\documentclass[11pt]{article}
\usepackage{graphicx}
\usepackage{amssymb}
\usepackage{amsmath}
\usepackage{setspace}

\newtheorem{theo}{Theorem}[section]
\newtheorem{lem}{Lemma}[section]
\newtheorem{prop}{Proposition}[section]
\newtheorem{cor}{Corollary}[section]
\newtheorem{rem}{Remark}[section]

\newenvironment{dem}[1][Proof]{\noindent {\itshape\ #1} }{\hfill $\square$}
\makeatletter
\def\@biblabel#1{\@ifnotempty{#1}{#1.}}
\input{tcilatex}

\begin{document}

\title{Lipschitz modulus of linear and convex systems with the Hausdorff
metric\thanks{%
This research has been partially supported by Grants
MTM2014-59179-C2-(1-2)-P and PGC2018-097960-B-C2(1,2) from MINECO/MICINN,
Spain, and FEDER "Una manera de hacer Europa", European Union.}}
\date{}
\author{G. Beer\thanks{%
Department of Mathematics, California State University Los Angeles, 5151
State University rive, Los Angeles, California 90032, USA
(gbeer@cslanet.calstatela.edu).}~~$\cdot $~~ M.J. C\'{a}novas\thanks{%
Center of Operations Research, Miguel Hern\'{a}ndez University of Elche,
03202 Elche (Alicante), Spain (canovas@umh.es, parra@umh.es).}~~$\cdot $%
~~M.A. L\'{o}pez\thanks{%
Department of Mathematics, University of Alicante, 03071 Alicante, Spain
(marco.antonio@ua.es).}~~$\cdot $~~J. Parra\footnotemark[3]  }
\maketitle

\begin{abstract}
\medskip This paper analyzes the Lipschitz behavior of the feasible set in
two parametric settings, associated with linear and convex systems in $%
\mathbb{R}^{n}.$ To start with, we deal with the parameter space of linear
(finite/semi-infinite) systems identified with the corresponding sets of
coefficient vectors, which are assumed to be closed subsets of $\mathbb{R}%
^{n+1}.$ In this framework, where the Hausdorff distance is used to measure
the size of perturbations, an explicit formula for computing the Lipschitz
modulus of the feasible set mapping is provided. As direct antecedent, we
appeal to its counterpart in the parameter space of all linear systems with
a fixed index set, $T,$ where the Chebyshev (pseudo) distance was considered
to measure the perturbations. Indeed, the stability (and, particularly,
Lipschitz properties) of linear systems in the Chebyshev framework has been
widely analyzed in the literature. Here, through an appropriate indexation
strategy, we take advantage of previous results to derive the new ones in
the Hausdorff setting. In a second stage, the possibility of perturbing
directly the set of coefficient vectors of a linear system allows us to
provide new contributions on the Lipschitz behavior of convex systems via
linearization techniques.

\noindent \textbf{Keywords. }Lipschitz modulus \textperiodcentered\ Feasible
set mapping \textperiodcentered\ Uncertain inequality systems
\textperiodcentered\ Hausdorff metric \textperiodcentered\ Indexation

\noindent \textbf{Mathematics Subject Classification: }90C31, 49J53, 49K40,
90C05, 90C25, 90C34.
\end{abstract}

\section{ Introduction}

This paper is initially focussed on the Lipschitz behavior of the \emph{%
feasible set }associated with a parametric family of \emph{linear inequality
systems }of the form:%
\begin{equation}
\{a^{\prime }x\leq b,\,\left( a,b\right) \in U\},U\in CL\left( \mathbb{R}%
^{n+1}\right) ,  \label{eq_problemRO}
\end{equation}%
where $x\in \mathbb{R}^{n}$ is the vector of variables, $CL\left( \mathbb{R}%
^{n+1}\right) $ is the \emph{parameter space }of all nonempty closed subsets
in $\mathbb{R}^{n+1}.$ Elements in $U\in CL\left( \mathbb{R}^{n+1}\right) $
are denoted as $\left( a,b\right) ,$ where $a\in \mathbb{R}^{n}$ and $b\in 
\mathbb{R}.$ Given $x,y\in \mathbb{R}^{n},$ $x^{\prime }y$ represents the
usual inner product of $x$ and $y.$ When $U$ is an infinite set, (\ref%
{eq_problemRO}) is a linear \emph{semi-infinite} inequality system. Observe
that, in this framework, perturbations fall on $U$ and, so, obviously, two
different systems, associated with different sets $U_{1},U_{2}\in CL\left( 
\mathbb{R}^{n+1}\right) ,$ can have different cardinality. This setting
includes as a particular case the parametric family of linear systems coming
from linearizing convex inequalities of the form%
\begin{equation}
f\left( x\right) \leq 0,\text{ }f\in \Gamma ,  \label{eq_convex_ineq}
\end{equation}%
where $\Gamma $ represents the set of all finite-valued convex functions on $%
\mathbb{R}^{n}.$ Specifically$,$ the feasible set of (\ref{eq_convex_ineq})
does coincide with the one of the linear system, 
\begin{equation}
\{a^{\prime }x\leq b:b=a^{\prime }z-f\left( z\right) ,\text{~~}\left(
z,a\right) \in \mathrm{gph}\partial f\},  \label{eq_linearization_system}
\end{equation}%
where $\mathrm{gph}\partial f$ represents the graph of the subdifferential
mapping $\partial f:\mathbb{R}^{n}\rightrightarrows \mathbb{R}^{n}$ given by 
\begin{equation*}
\partial f\left( z\right) :=\{a\in \mathbb{R}^{n}\mid f\left( x\right) \geq
f\left( z\right) +a^{\prime }\left( x-z\right) ,\text{ for all }x\in \mathbb{%
R}^{n}\}.
\end{equation*}%
\bigskip

The main objectives of this work consist of analyzing the Lipschitzian
behavior of the parametrized linear system (\ref{eq_problemRO}) and to apply
the obtained results to derive new contributions on the convex case (\ref%
{eq_convex_ineq}) via the standard linearization (\ref%
{eq_linearization_system}). We emphasize the fact that previous results
about stability of subdifferential (traced out from \cite{BCLP19}) are also
used in the study of this convex case.

Formally, associated with (\ref{eq_problemRO}), we consider the \emph{%
feasible set mapping}, $\mathcal{F}:CL\left( \mathbb{R}^{n+1}\right)
\rightrightarrows \mathbb{R}^{n},$ which assigns to each $U\in CL\left( 
\mathbb{R}^{n+1}\right) $ the set of solutions of the corresponding system: 
\begin{equation}
\mathcal{F}\left( U\right) :=\left\{ x\in \mathbb{R}^{n}:a^{\prime }x\leq b%
\text{ for all }\left( a,b\right) \in U\right\} ,\text{ }U\in CL\left( 
\mathbb{R}^{n+1}\right) .  \label{eq_0004}
\end{equation}%
The parameter space, $CL\left( \mathbb{R}^{n+1}\right) ,$ will be endowed
with Hausdorff distance (see Section 2 for details). For convenience, we
deal with closed sets, but the study could be carried out with general
nonempty sets, since both the feasible set mapping and the Hausdorff
distance (pseudo-distance in such a case) do not distinguish between
nonempty sets and their closures.

The main original contributions of the present paper consists of providing a
formula for computing the \emph{Lipschitz modulus} of $\mathcal{F}$ at $%
\left( U_{0},x_{0}\right) \in \mathrm{gph}\mathcal{F}$ and, in a second
stage, to derive a Lipschitzian type condition for the feasible set of the
parametrized convex system (\ref{eq_convex_ineq}). Roughly speaking we
provide measures (or estimations) of the rate of variation of feasible
points, around a nominal one $x_{0}\in \mathbb{R}^{n}$ with respect to
perturbations of the nominal parameter ($U_{0}\in CL\left( \mathbb{R}%
^{n+1}\right) $ in the case of linear systems and $f_{0}\in \Gamma $ in the
convex case).

As immediate antecedents of the present work we cite \cite{CGP08} (see also
updated results in \cite{CGHP19}) and \cite{CHLP18}. The first paper
computes the Lipschitz modulus of the feasible set mapping in the context of
linear systems with an arbitrarily fixed index set $T$ of the form%
\begin{equation}
\left\{ a_{t}^{\prime }x\leq b_{t},\,t\in T\right\} ,  \label{eq_system_T}
\end{equation}%
where $x\in \mathbb{R}^{n}$ is the variable and $\left( a_{t},b_{t}\right)
_{t\in T}\in \left( \mathbb{R}^{n+1}\right) ^{T}.$ The parameter space
considered there, $\left( \mathbb{R}^{n+1}\right) ^{T},$ is formed by all
functions from $T$ to $\mathbb{R}^{n+1}$ and it is endowed with the
(extended) Chebyshev distance. The reader is addressed to the monograph \cite%
[Chapter 6]{libro} for a comprehensive study of such systems. The results of 
\cite{CGP08} do not apply directly to our current setting unless some
appropriate connection between both parameter spaces, $CL\left( \mathbb{R}%
^{n+1}\right) $ and $\left( \mathbb{R}^{n+1}\right) ^{T},$ was established.
In relation to this point, we appeal to paper \cite{CHLP18}, which provides
the motivation and background from the methodological point of view. That
paper is focussed on the \emph{calmness modulus }(see again Section 2), and
takes advantage of previous results developed in the context of systems (\ref%
{eq_system_T}), to derive new contributions for the parametrized system (\ref%
{eq_problemRO}). Formally, \cite{CHLP18} introduces an appropriate \emph{%
indexation mapping, }assigning to each set in $CL\left( \mathbb{R}%
^{n+1}\right) $ an element in $\left( \mathbb{R}^{n+1}\right) ^{T}$ in such
a way that the Hausdorff distance around $U_{0}\in CL\left( \mathbb{R}%
^{n+1}\right) $ translates into the Chebyshev distance around its image in $%
CL\left( \mathbb{R}^{n+1}\right) $. That indexation strategy is shown to be
inappropriate for studying the Lipschitz (instead of calmness) modulus and,
in relation to this fact, a new indexation strategy is introduced in Section
3.

The problem of analyzing the relationship among different parametric
contexts was also addressed in \cite{CLP05} and \cite{CLP08} from a
different perspective, mainly focussed on the lower semicontinuity of the
feasible set mapping.

Now we summarize the structure of the paper. Section 2 gathers some
definitions and key results of the background on the Lipschitz modulus in
the context of systems (\ref{eq_system_T}), indexations, and stability of
subdifferentials. Section 3 develops the study of the Lipschitz modulus of $%
\mathcal{F},$ including the definition of an appropriate indexation which
allows us to take advantage of the background about systems (\ref%
{eq_system_T}). Finally, Section 4 applies the results of previous section
to tackle the convex case.

\section{Preliminaries and first results}

To start with, recall that a set-valued mapping $\mathcal{M}%
:Y\rightrightarrows X$ between metric spaces (both distances denoted by $d$)
has the \emph{Aubin property} (also called pseudo-Lipschitz --cf. \cite%
{KlKu02}-- or Lipschitz-like --cf. \cite{mor06a}--) at $\left(
y_{0},x_{0}\right) \in \mathrm{gph}\mathcal{M}$ if there exist a constant $%
\kappa \geq 0$ and neighborhoods $W$ of $x_{0}$ and $V$ of $y_{0}$ such that 
\begin{equation}
d\left( x_{1},\mathcal{M}\left( y_{2}\right) \right) \leq \kappa d\left(
y_{1},y_{2}\right) \text{ for all }y_{1},y_{2}\subset V\text{ and all }%
x_{1}\in \mathcal{M}\left( y_{1}\right) \cap W.  \label{eq_Aubin}
\end{equation}%
The infimum of constants $\kappa $ over all $\left( \kappa ,W,V\right) $
satisfying (\ref{eq_Aubin}) is called the \emph{Lipschitz modulus} of $%
\mathcal{M}$ at $\left( y_{0},x_{0}\right) ,$ denoted by $\mathrm{lip}%
\mathcal{M}\left( y_{0},x_{0}\right) ,$ and it is defined as $+\infty $ when
the Aubin property fails at $\left( y_{0},x_{0}\right) .$ The Aubin property
of $\mathcal{M}$ at $\left( y_{0},x_{0}\right) \in \mathrm{gph}\mathcal{M}$
is known to be equivalent to the metric regularity of its inverse mapping $%
\mathcal{M}^{-1}$ at $\left( x_{0},y_{0}\right) \in \mathrm{gph}\mathcal{M};$
moreover, $\mathrm{lip}\mathcal{M}\left( y_{0},x_{0}\right) $ is known to
coincide with the modulus of metric regularity of $\mathcal{M}^{-1}$ at $%
\left( x_{0},y_{0}\right) .\,$So, we can write 
\begin{equation*}
\mathrm{lip}\mathcal{M}\left( y_{0},x_{0}\right) =\underset{\left(
x,y\right) \rightarrow \left( x_{0},y_{0}\right) }{\limsup }\frac{d\left( x,%
\mathcal{M}(y)\right) }{d\left( y,\mathcal{M}^{-1}(x)\right) },
\end{equation*}%
under the conventions $\frac{0}{0}:=0,$ $\frac{1}{0}:=+\infty ,$ and $\frac{1%
}{+\infty }:=0.$

The particularization of (\ref{eq_Aubin}) to $y_{2}=y_{0}$ yields the
definition of \emph{calmness }of $\mathcal{M}$ at $\left( y_{0},x_{0}\right)
,$ whose associated \emph{calmness modulus}, $\mathrm{clm}\mathcal{M}\left(
y_{0},x_{0}\right) ,$ is defined analogously. It is also known that the
calmness of $\mathcal{M}$ at $\left( y_{0},x_{0}\right) $ is equivalent to
the metric subregularity of $\mathcal{M}^{-1}$ at $\left( x_{0},y_{0}\right) 
$, and that the corresponding moduli do coincide; so, 
\begin{equation*}
\mathrm{clm}\mathcal{M}\left( y_{0},x_{0}\right) =\underset{x\rightarrow
x_{0}}{\limsup }\frac{d\left( x,\mathcal{M}(y_{0})\right) }{d\left( y_{0},%
\mathcal{M}^{-1}(x)\right) }.
\end{equation*}%
Clearly $\mathrm{clm}\mathcal{M}\left( y^{0},x_{0}\right) \leq \mathrm{lip}%
\mathcal{M}\left( y^{0},x_{0}\right) .$ For additional information about the
Aubin property, calmness, and related topics of variational analysis, the
reader is addressed to \cite{DoRo, KlKu02, mor06a, rw}.

\subsection{Indexation strategies and calmness of linear systems}

For comparative purposes and as a motivation of the results of Section 3,
this subsection recalls some details about the indexation introduced in \cite%
{CHLP18}. First, we fix the topologies considered in the space of variables, 
$\mathbb{R}^{n},$ and the parameter spaces, $CL\left( \mathbb{R}%
^{n+1}\right) $ and $\left( \mathbb{R}^{n+1}\right) ^{T}$.

Unless otherwise stated, $\mathbb{R}^{n}$ (space of variables) is equipped
with an arbitrary norm, $\left\Vert \cdot \right\Vert ,$ while $\mathbb{R}%
^{n+1}$ (space of coefficient vectors of linear systems) is endowed with the
norm 
\begin{equation}
\left\Vert \left( a,b\right) \right\Vert =\max \left\{ \left\Vert
a\right\Vert _{\ast },\left\vert b\right\vert \right\} ,\text{ }\left(
a,b\right) \in \mathbb{R}^{n+1},  \label{eq_normRn+1}
\end{equation}%
where $\left\Vert \cdot \right\Vert _{\ast }$ represents the dual norm of $%
\left\Vert \cdot \right\Vert $ in $\mathbb{R}^{n}$, which is given by%
\begin{equation*}
\left\Vert a\right\Vert _{\ast }=\sup_{\left\Vert x\right\Vert \leq
1}a^{\prime }x.
\end{equation*}

The space $CL\left( \mathbb{R}^{n+1}\right) $ is endowed with the (extended)
Hausdorff distance $d_{H}:CL\left( \mathbb{R}^{n+1}\right) \times CL\left( 
\mathbb{R}^{n+1}\right) \rightarrow \lbrack 0,+\infty ]$ given by 
\begin{equation*}
d_{H}\left( U_{1},U_{2}\right) :=\max \{e\left( U_{1},U_{2}\right) ,e\left(
U_{2},U_{1}\right) \},
\end{equation*}%
where $e\left( U_{i},U_{j}\right) ,$ $i,j=1,2,$ represents the excess of $%
U_{i}$ over $U_{j},$ 
\begin{equation}
e\left( U_{i},U_{j}\right) :=\inf \left\{ \varepsilon >0:U_{i}\subset
U_{j}+\varepsilon \mathbb{B}\right\} =\sup \left\{ d\left( x,U_{j}\right)
:x\in U_{i}\right\} ,  \notag
\end{equation}%
where $\mathbb{B}$ denotes the closed unit ball in $\mathbb{R}^{n+1}.$ See 
\cite[Section 3.2]{Beer} for details about the Hausdorff distance in general
settings.

In $\left( \mathbb{R}^{n+1}\right) ^{T},$ the (extended) Chebyshev (or
supremum) distance, $d_{\infty }:\left( \mathbb{R}^{n+1}\right) ^{T}\times
\left( \mathbb{R}^{n+1}\right) ^{T}\rightarrow \lbrack 0,+\infty ],$ given
by 
\begin{equation}
d_{\infty }\left( \sigma _{1},\sigma _{2}\right) :=\sup_{t\in T}\left\Vert
\sigma _{1}\left( t\right) -\sigma _{2}\left( t\right) \right\Vert ,
\label{eq_def_d_inf}
\end{equation}%
is considered.

As commented in the introduction, paper \cite{CHLP18} analyzes the calmness
of $\mathcal{F}$ at $\left( U_{0},x_{0}\right) \in \mathrm{gph}\mathcal{F}$
via the calmness of the feasible set mapping associated with systems (\ref%
{eq_system_T}), $\mathcal{F}^{T}:\left( \mathbb{R}^{n+1}\right)
^{T}\rightrightarrows \mathbb{R}^{n},$ which is given by 
\begin{equation}
\mathcal{F}^{T}\left( \sigma \right) :=\left\{ x\in \mathbb{R}%
^{n}:a_{t}^{\prime }x\leq b_{t},\text{ }t\in T\right\} ,\text{ }\sigma
=\left( a_{t},b_{t}\right) _{t\in T}\in \left( \mathbb{R}^{n+1}\right) ^{T}.
\label{eq_0005}
\end{equation}%
To do this, a particular \emph{indexation mapping }between $CL\left( \mathbb{%
R}^{n+1}\right) $ and $\left( \mathbb{R}^{n+1}\right) ^{T}$ is introduced.
Recall that $\sigma \in \left( \mathbb{R}^{n+1}\right) ^{T}$ is said to be
an \emph{indexation} of $U\in CL\left( \mathbb{R}^{n+1}\right) $ if 
\begin{equation*}
\mathrm{rge}\left( \sigma \right) =U,
\end{equation*}%
where $\mathrm{\,`rge}$' means range (or image); specifically, \cite{CHLP18}
considers $T:=\mathbb{R}^{n+1}\ $and assigns to each $U\in CL\left( \mathbb{R%
}^{n+1}\right) $ an indexation $\mathcal{I}_{U}\in \left( \mathbb{R}%
^{n+1}\right) ^{\mathbb{R}^{n+1}}$ defined as 
\begin{equation}
\mathcal{I}_{U}\left( t\right) :=\left\{ 
\begin{tabular}{cc}
$t$ & $\text{if }t\in U,$ \\ 
$P_{U}\circ P_{U_{0}}(t)$ & $\text{if }t\notin U,$%
\end{tabular}%
\right.  \label{eq_indexation}
\end{equation}%
where, for each $U_{1}\in $ $CL\left( \mathbb{R}^{n+1}\right) ,$ $P_{U_{1}}:%
\mathbb{R}^{n+1}\rightarrow \mathbb{R}^{n+1}$ is a particular selection of
the \emph{metric projection} multifunction on $U_{1};$ i.e., $%
P_{U_{1}}\left( t\right) $ is a best approximation of $t\in \mathbb{R}^{n+1}$
on $U_{1}\in CL\left( \mathbb{R}^{n+1}\right) .$ Observe that, in
particular, $\mathcal{I}_{U_{0}}=P_{U_{0}}.$ A comparative analysis with
other possible indexations, and particularly one given in \cite{ChanMar}, is
carried out in \cite[Section 3]{CHLP18}. Theorem 3.1 in \cite{CHLP18} shows
that 
\begin{equation}
d_{\infty }\left( \mathcal{I}_{U},\mathcal{I}_{U_{0}}\right) =d_{H}\left(
U,U_{0}\right) \text{ for all }U\in CL\left( \mathbb{R}^{n+1}\right) .
\label{eq_0003}
\end{equation}%
Example 3.1 in the same paper shows that $U\mapsto P_{U}$ is not an adequate
indexation mapping in relation to calmness, as far as Chebyshev distances
between projections, $d_{\infty }\left( P_{U},P_{U_{0}}\right) ,$ can be
much larger than Hausdorff distances between sets $d_{H}\left(
U,U_{0}\right) .$

Indexation mapping $\mathcal{I}$ in (\ref{eq_indexation}) is suitable for
the study of the calmness property of $\mathcal{F},$ but it is no longer
enough for the Aubin property, for which we need more, namely: $d_{\infty
}\left( \sigma _{1},\sigma _{2}\right) =d_{H}\left( U_{1},U_{2}\right) $
when $\sigma _{1}$ and $\sigma _{2}$ are indexations of two sets $U_{1}$ and 
$U_{2}$ close enough to the nominal set $U_{0}.$ The price to pay is that
the definition of $\sigma _{2}$ depends not only on $U_{0},$ but also on $%
U_{1}.$ Such an indexation strategy is defined in the proof of Theorem \ref%
{theo_lip_U} and constitutes one of the main contributions of the this
section.

\subsection{On the stability of subdifferentials}

This subsection gathers some stuff traced out from \cite{BCLP19} about
stability of subdifferentials of convex functions at a point $x_{0}\in 
\mathbb{R}^{n}$ and provides some extensions and consequences on the
stability over a compact set $K_{0}\subset \mathbb{R}^{n}.$ This results
will be used in Section 4.

Given any two functions $f_{1},f_{2}\in \Gamma $ and a compact subset $%
K\subset \mathbb{R}^{n}$ we use the notation 
\begin{equation}
d_{K}\left( f_{1},f_{2}\right) =\sup_{x\in K}\left\vert f_{1}\left( x\right)
-f_{2}\left( x\right) \right\vert .  \label{eq_0007}
\end{equation}%
The following theorem gathers two stability conditions for subdifferentials.
The first one, which is a direct consequence of \cite[Theorem 24.5]{Rock},
provides the Hausdorff upper semicontinuity of the multifunction which
assigns to each pair $\left( f,x\right) \in \Gamma \times \mathbb{R}^{n}$
the subdifferential of $f$ at $x,$ 
\begin{equation*}
\partial f\left( x\right) :=\{u\in \mathbb{R}^{n}\mid f\left( z\right) \geq
f\left( x\right) +u^{\prime }\left( z-x\right) ,\text{ for all }z\in \mathbb{%
R}^{n}\}.
\end{equation*}%
On the other hand, condition $\left( ii\right) $ expresses a certain uniform
lower H\"{o}lder type property.

\begin{theo}
\label{Th_lower}\label{Th_Subdif}Let $x_{0}\in \mathbb{R}^{n},$ $\alpha >0,$
and $K:=x_{0}+\alpha \mathbb{B}.$ One has:

$\left( i\right) $ \emph{\cite[Prop. 2.1]{BCLP19} }Given $f_{0}\in \Gamma $
and $\varepsilon >0,$ there exists $\delta >0$ such that 
\begin{equation*}
\partial f\left( x_{0}+\delta \mathbb{B}\right) \subset \partial f_{0}\left(
x_{0}\right) +\varepsilon \mathbb{B},
\end{equation*}%
provided that $f\in \Gamma $ satisfies $d_{K}\left( f,f_{0}\right) \leq
\delta .\medskip $

$\left( ii\right) $\emph{\cite[Thm. 3.4]{BCLP19} }For any $0<\delta \leq
\alpha ^{2},$ and any $f_{1},f_{2}\in \Gamma $ such that $%
d_{K}(f_{1},f_{2})\leq \delta ,$ we have 
\begin{equation}
\partial f_{1}\left( x_{0}\right) \subset \partial f_{2}\left( x_{0}+\sqrt{%
\delta }\mathbb{B}\right) +4\sqrt{\delta }\mathbb{B}.  \label{marco101}
\end{equation}
\end{theo}

\begin{cor}
\label{Cor_subdifferential}Let $K_{0}\subset \mathbb{R}^{n}$ a compact set, $%
\alpha >0,$ and $K:=K_{0}+\alpha \mathbb{B}.$ One has:

$\left( i\right) $ Given $f_{0}\in \Gamma $ and $\varepsilon >0,$ there
exists $\delta >0$ such that 
\begin{equation*}
\partial f\left( K_{0}+\delta \mathbb{B}\right) \subset \partial f_{0}\left(
K_{0}\right) +\varepsilon \mathbb{B},
\end{equation*}%
provided that $f\in \Gamma $ satisfies $d_{K}\left( f,f_{0}\right) \leq
\delta .$

$\left( ii\right) $ For any $0<\delta \leq \alpha ^{2},$ and any $%
f_{1},f_{2}\in \Gamma $ such that $d_{K}(f_{1},f_{2})\leq \delta ,$ we have 
\begin{equation*}
\partial f_{1}\left( K_{0}\right) \subset \partial f_{2}\left( K_{0}+\sqrt{%
\delta }\mathbb{B}\right) +4\sqrt{\delta }\mathbb{B}.
\end{equation*}

$\left( iii\right) $ Given $f_{0}\in \Gamma $ and $\varepsilon >0,$ there
exists $\delta _{0}>0$ such that for any $0<\delta \leq \delta _{0},$ and
any $f\in \Gamma ,$ with $d_{K}\left( f,f_{0}\right) \leq \delta ,$ one has%
\begin{equation*}
d_{H}\left( \partial f\left( K_{0}+\sqrt{\delta }\mathbb{B}\right) ,\partial
f_{0}\left( K_{0}\right) \right) \leq \varepsilon .
\end{equation*}
\end{cor}

\begin{dem}
$\left( i\right) $ follows the same argument of the proof of Theorem \ref%
{Th_lower}$\left( i\right) .$ Here we present a sketch for completeness.
Arguing by contradiction, assume the existence of sequences $\left\{
f_{r}\right\} \subset \Gamma $ and $\{\left( x_{r},u_{r}\right) \}\subset 
\mathbb{R}^{n}\times \mathbb{R}^{n}$ such that $d_{K}\left(
f_{r},f_{0}\right) \leq \frac{1}{r},$ $x_{r}\in K_{0}+\frac{1}{r}\mathbb{B},$
and $u_{r}\in \partial f_{r}\left( x_{r}\right) \diagdown \partial
f_{0}\left( x_{0}\right) +\varepsilon \mathbb{B},$ $r=1,2,...$ In this way
we reach a contradiction (see \cite[Theorem 24.5]{Rock}) as far as we may
assume without loss of generality that $\{x_{r}\}$ converges to a certain $%
x_{0}\in K_{0}.$

$\left( ii\right) $ Comes straightforwardly from \ref{Th_lower}$\left(
ii\right) .$ Indeed, let $0<\delta \leq \alpha ^{2},$ let $f_{1},f_{2}\in
\Gamma $ be such that $d_{K}(f_{1},f_{2})\leq \delta ,$ and take any $%
x_{0}\in K_{0}.\,\ $We have $d_{x_{0}+\alpha \mathbb{B}}(f_{1},f_{2})\leq
d_{K}(f_{1},f_{2})\leq \delta ,$ which entails 
\begin{equation*}
\partial f_{1}\left( x_{0}\right) \subset \partial f_{2}\left( x_{0}+\sqrt{%
\delta }\mathbb{B}\right) +4\sqrt{\delta }\mathbb{B\subset }\partial
f_{2}\left( K_{0}+\sqrt{\delta }\mathbb{B}\right) +4\sqrt{\delta }\mathbb{B}.
\end{equation*}

$\left( iii\right) $ Take $f_{0}\in \Gamma $ and $\varepsilon >0.$ From
condition $\left( i\right) \ $there exists $\delta _{1}>0$ such that 
\begin{equation*}
\partial f\left( K_{0}+\delta _{1}\mathbb{B}\right) \subset \partial
f_{0}\left( K_{0}\right) +\varepsilon \mathbb{B},\text{ provided that }%
d_{K}\left( f,f_{0}\right) \leq \delta _{1},\text{ }f\in \Gamma .
\end{equation*}%
We may assume $\delta _{1}\leq 1.$

Define%
\begin{equation*}
0<\delta _{0}:=\min \{\delta _{1}^{2},\alpha ^{2},\left( \frac{\varepsilon }{%
4}\right) ^{2}\},
\end{equation*}%
and take $0<\delta \leq \delta _{0},$ and $f\in \Gamma $ such that $%
d_{K}\left( f,f_{0}\right) \leq \delta .$ Then, since $\sqrt{\delta }\leq 
\sqrt{\delta _{0}}\leq \delta _{1}$ and $\delta \leq \delta _{1}^{2}\leq
\delta _{1}$ (which yields $d_{K}\left( f,f_{0}\right) \leq \delta _{1}),$
one has$,$ 
\begin{equation*}
\partial f\left( K_{0}+\sqrt{\delta }\mathbb{B}\right) \subset \partial
f\left( K_{0}+\delta _{1}\mathbb{B}\right) \subset \partial f_{0}\left(
K_{0}\right) +\varepsilon \mathbb{B}.
\end{equation*}%
On the other hand since $0<\delta \leq \alpha ^{2},$ condition $\left(
ii\right) $ yields 
\begin{equation*}
\partial f_{0}\left( K_{0}\right) \subset \partial f\left( K_{0}+\sqrt{%
\delta }\mathbb{B}\right) +4\sqrt{\delta }\mathbb{B\subset }\partial f\left(
K_{0}+\sqrt{\delta }\mathbb{B}\right) +\varepsilon \mathbb{B},
\end{equation*}%
where the last inclusion comes from $\delta \leq \left( \frac{\varepsilon }{4%
}\right) ^{2}.$
\end{dem}

\section{Lipschitz modulus of $\mathcal{F}$ in the Hausdorff setting}

This section provides a point-based formula for $\mathrm{lip\,}\mathcal{F}%
\left( U_{0},x_{0}\right) $ through a previously established expression of $%
\mathrm{lip\,}\mathcal{F}^{T}\left( \sigma _{0},x_{0}\right) ,$ which is
recalled in the next theorem$.$ First, we introduce some notation: Given $%
X\subset \mathbb{R}^{k},$ $k\in \mathbb{N},$ we denote by \textrm{conv}$X$
the \emph{convex hull }of $X,$ and $\mathrm{int}X,$ $\mathrm{cl}X$ and $%
\mathrm{bd}X$ stand, respectively, for the interior, the closure and the
boundary of $X.$

\begin{theo}
\emph{(see \cite[Theorem 1$\left( ii\right) $]{CGP08})}\label{TH_lip} Let $%
\left( \sigma _{0},x_{0}\right) \in \mathrm{gph}\mathcal{F}^{T},$ with $%
\sigma _{0}=\left( a_{t}^{0},b_{t}^{0}\right) _{t\in T}\in \left( \mathbb{R}%
^{n+1}\right) ^{T}.$ Assume that $\left\{ a_{t}^{0},\text{ }t\in T\right\} $
is bounded. Then 
\begin{equation*}
\mathrm{lip\,}\mathcal{F}^{T}\left( \sigma _{0},x_{0}\right) =\frac{%
\left\Vert x_{0}\right\Vert +1}{d_{\ast }\left( 0_{n},C_{0}\right) },
\end{equation*}%
where 
\begin{equation*}
C_{0}:=\left\{ u\in \mathbb{R}^{n}:\left( u,u^{\prime }x_{0}\right) \in 
\mathrm{cl\,conv\,}\left\{ \left( a_{t}^{0},b_{t}^{0}\right) ,~t\in
T\right\} \right\} .
\end{equation*}
\end{theo}

The following lemma constitutes a key step for deriving the announced
formula for $\mathrm{lip}\mathcal{F}\left( U_{0},x_{0}\right) $. In it we
construct appropriate indexations of sets $U_{1},U_{2}\in CL\left( \mathbb{R}%
^{n+1}\right) ,$ denoted by $\sigma _{1},\sigma _{2},$ which preserve the
distance between them; i.e., $d_{H}\left( U_{1},U_{2}\right) =d_{\infty
}\left( \sigma _{1},\sigma _{2}\right) ,$ and and we obtain Lipschitz
estimates for each $d_{\infty }\left( \sigma _{i},\sigma _{0}\right) $ in
terms of $d_{H}\left( U_{j},U_{0}\right) ,$ $j=1,2.$ At this moment, we
observe that, in general, given any\emph{\ }$\sigma _{1},\sigma _{2}\in
\left( \mathbb{R}^{n+1}\right) ^{T},$ $T$ being arbitrary, one has\emph{\ }%
\begin{equation}
d_{H}\left( \mathrm{cl\,rge\,}\sigma _{1},\mathrm{cl\,rge\,}\sigma
_{2}\right) \leq d_{\infty }\left( \sigma _{1},\sigma _{2}\right) .
\label{q_0006}
\end{equation}%
For simplicity in the notation, in the lemma let us write $P_{i}\left(
t\right) $ instead $P_{U_{i}}\left( t\right) ,$ for $t\in \mathbb{R}^{n+1}$
and $i=0,1,2.$

\begin{lem}
Let $U_{0}\in CL\left( \mathbb{R}^{n+1}\right) ,$ $T:=\mathbb{R}^{n+1},$ and 
$\sigma _{0}:=P_{0}\in \left( \mathbb{R}^{n+1}\right) ^{T}.$ Associated with
each pair of subsets $\left( U_{1},U_{2}\right) \in CL\left( \mathbb{R}%
^{n+1}\right) \times CL\left( \mathbb{R}^{n+1}\right) $ let us define a pair
of functions $\left( \sigma _{1},\sigma _{2}\right) \in \left( \mathbb{R}%
^{n+1}\right) ^{T}\times \left( \mathbb{R}^{n+1}\right) ^{T}$ as follows:
for each $t\in \mathbb{R}^{n+1}$, 
\begin{equation*}
\left( \sigma _{1},\sigma _{2}\right) \left( t\right) :=\left\{ 
\begin{tabular}{l}
$\left( P_{1}\left( t\right) ,P_{2}\left( t\right) \right) ,$ if $t\in
U_{1}\cup U_{2},$ \\ 
$\left( P_{1}\circ P_{0}\left( t\right) ,\text{ }P_{2}\circ P_{1}\circ
P_{0}\left( t\right) \right) ,$ if $t\notin U_{1}\cup U_{2}.$%
\end{tabular}%
\medskip \right.
\end{equation*}%
Then we have%
\begin{eqnarray*}
d_{\infty }\left( \sigma _{i},\sigma _{0}\right) &\leq &3\max \{d_{H}\left(
U_{1},U_{0}\right) ,d_{H}\left( U_{2},U_{0}\right) \},\text{ }i=1,2, \\
d_{\infty }\left( \sigma _{1},\sigma _{2}\right) &=&d_{H}\left(
U_{1},U_{2}\right) .
\end{eqnarray*}
\end{lem}

\begin{dem}
Take any $\left( U_{1},U_{2}\right) \in CL\left( \mathbb{R}^{n+1}\right)
\times CL\left( \mathbb{R}^{n+1}\right) ,$ and the associated $\left( \sigma
_{1},\sigma _{2}\right) \in \left( \mathbb{R}^{n+1}\right) ^{T}\times \left( 
\mathbb{R}^{n+1}\right) ^{T}$ as in the statement of the lemma. First, let
us see that 
\begin{equation}
d_{\infty }\left( \sigma _{1},\sigma _{0}\right) \leq 3\max \{d_{H}\left(
U_{1},U_{0}\right) ,d_{H}\left( U_{2},U_{0}\right) \}.  \label{e_001}
\end{equation}%
For $t\in U_{1}$ we have 
\begin{equation*}
\left\Vert \sigma _{1}\left( t\right) -\sigma _{0}\left( t\right)
\right\Vert =\left\Vert t-P_{0}\left( t\right) \right\Vert \leq e\left(
U_{1},U_{0}\right) \leq d_{H}\left( U_{1},U_{0}\right) .
\end{equation*}%
For $t\in U_{2}$ we have 
\begin{eqnarray*}
\left\Vert \sigma _{1}\left( t\right) -\sigma _{0}\left( t\right)
\right\Vert &\leq &\left\Vert \sigma _{1}\left( t\right) -\sigma _{2}\left(
t\right) \right\Vert +\left\Vert \sigma _{2}\left( t\right) -\sigma
_{0}\left( t\right) \right\Vert \\
&=&\left\Vert P_{1}\left( t\right) -t\right\Vert +\left\Vert t-P_{0}\left(
t\right) \right\Vert \\
&\leq &e\left( U_{2},U_{1}\right) +e\left( U_{2},U_{0}\right) \leq
2d_{H}\left( U_{2},U_{0}\right) +d_{H}\left( U_{1},U_{0}\right) \\
&\leq &3\max \{d_{H}\left( U_{1},U_{0}\right) ,d_{H}\left(
U_{2},U_{0}\right) \}.
\end{eqnarray*}%
For $t\notin U_{1}\cup U_{2}$ we have 
\begin{equation*}
\left\Vert \sigma _{1}\left( t\right) -\sigma _{0}\left( t\right)
\right\Vert =\left\Vert P_{1}\circ P_{0}\left( t\right) -P_{0}\left(
t\right) \right\Vert \leq e\left( U_{0},U_{1}\right) \leq d_{H}\left(
U_{1},U_{0}\right) .
\end{equation*}%
In summary, (\ref{e_001}) holds in any case.

Now, let us check that 
\begin{equation}
d_{\infty }\left( \sigma _{2},\sigma _{0}\right) \leq 3\max \{d_{H}\left(
U_{1},U_{0}\right) ,d_{H}\left( U_{2},U_{0}\right) \}.  \label{e_002}
\end{equation}%
For $t\in U_{1}\cup U_{2}$ the arguments are completely analogous to those
of $\sigma _{1}.$ For $t\notin U_{1}\cup U_{2}$ we have 
\begin{eqnarray*}
\left\Vert \sigma _{2}\left( t\right) -\sigma _{0}\left( t\right)
\right\Vert &\leq &\left\Vert \sigma _{2}\left( t\right) -\sigma _{1}\left(
t\right) \right\Vert +\left\Vert \sigma _{1}\left( t\right) -\sigma
_{0}\left( t\right) \right\Vert \\
&=&\left\Vert P_{2}\circ P_{1}\circ P_{0}\left( t\right) -P_{1}\circ
P_{0}\left( t\right) \right\Vert +\left\Vert \sigma _{1}\left( t\right)
-\sigma _{0}\left( t\right) \right\Vert \\
&\leq &e\left( U_{1},U_{2}\right) +d_{H}\left( U_{1},U_{0}\right) \leq
2d_{H}\left( U_{1},U_{0}\right) +d_{H}\left( U_{2},U_{0}\right) \\
&\leq &3\max \{d_{H}\left( U_{1},U_{0}\right) ,d_{H}\left(
U_{2},U_{0}\right) \}.
\end{eqnarray*}

So, we have established (\ref{e_002}).

The last step consists of checking $d_{\infty }\left( \sigma _{1},\sigma
_{2}\right) =d_{H}\left( U_{1},U_{2}\right) .$ On the one hand, 
\begin{equation*}
\sup_{t\in U_{1}}\left\Vert \sigma _{1}\left( t\right) -\sigma _{2}\left(
t\right) \right\Vert =\sup_{t\in U_{1}}\left\Vert t-P_{2}\left( t\right)
\right\Vert =e\left( U_{1},U_{2}\right) ,
\end{equation*}%
and, analogously, $\sup_{t\in U_{2}}\left\Vert \sigma _{1}\left( t\right)
-\sigma _{2}\left( t\right) \right\Vert =e\left( U_{2},U_{1}\right) .$ On
the other hand, for all $t\notin U_{1}\cup U_{2}$ we have 
\begin{equation*}
\left\Vert \sigma _{1}\left( t\right) -\sigma _{2}\left( t\right)
\right\Vert =\left\Vert P_{1}\circ P_{0}\left( t\right) -P_{2}\circ
P_{1}\circ P_{0}\left( t\right) \right\Vert \leq e\left( U_{1},U_{2}\right) .
\end{equation*}
\end{dem}

\begin{theo}
\label{theo_lip_U}Let $U_{0}\in CL\left( \mathbb{R}^{n+1}\right) $ be such
that $\left\{ a\in \mathbb{R}^{n}:\left( a,b\right) \in U_{0}\right\} $ is
bounded. Then, 
\begin{equation}
\mathrm{lip\,}\mathcal{F}\left( U_{0},x_{0}\right) =\frac{\left\Vert
x_{0}\right\Vert +1}{d_{\ast }\left( 0_{n},C_{U_{0}}\right) },  \label{eq_8}
\end{equation}%
where 
\begin{equation*}
C_{U_{0}}:=\left\{ u\in \mathbb{R}^{n}:\left( u,u^{\prime }x_{0}\right) \in 
\mathrm{cl\,conv\,}U_{0}\right\} .
\end{equation*}
\end{theo}

\begin{dem}
For simplicity, let us denote by $\kappa _{0}$ the right-hand side of (\ref%
{eq_8}). Take the indexation of $U_{0},$ $\sigma _{0}:=P_{U_{0}},$ and
observe that the pair $\left( \sigma _{0},x_{0}\right) $ satisfies all the
hypotheses of Theorem \ref{TH_lip}. Accordingly, $\mathrm{lip\,}\mathcal{F}%
^{T}\left( \sigma _{0},x_{0}\right) =\kappa _{0}.$ The rest of the proof is
devoted to showing that $\mathrm{lip\,}\mathcal{F}^{T}\left( \sigma
_{0},x_{0}\right) =\mathrm{lip\,}\mathcal{F}\left( U_{0},x_{0}\right) .$ Let 
$\varepsilon >0$ be arbitrarily given. By the definition of $\mathrm{lip\,}%
\mathcal{F}^{T}\left( \sigma _{0},x_{0}\right) $, there exists $\delta >0$
such that, for all $\sigma _{1},\sigma _{2}\in \left( \mathbb{R}%
^{n+1}\right) ^{T}$ with $d_{\infty }\left( \sigma _{i},\sigma _{0}\right)
<\delta ,$ $i=1,2,$ and all $x^{1}\in \mathcal{F}^{T}\left( \sigma
_{1}\right) $ with $\left\Vert x^{1}-x_{0}\right\Vert <\delta ,$ one has 
\begin{equation}
d\left( x^{1},\mathcal{F}^{T}\left( \sigma _{2}\right) \right) \leq \left(
\kappa _{0}+\varepsilon \right) d_{\infty }\left( \sigma _{1},\sigma
_{2}\right) .  \label{eq_12}
\end{equation}%
We are going to prove that, for all $U_{1},U_{2}\in CL\left( \mathbb{R}%
^{n+1}\right) $ with $d_{H}\left( U_{i},U_{0}\right) <\delta /3,$ $i=1,2,$
and all $x^{1}\in \mathcal{F}\left( U_{1}\right) $ with $\left\Vert
x^{1}-x_{0}\right\Vert <\delta ,$ one has 
\begin{equation}
d\left( x^{1},\mathcal{F}\left( U_{2}\right) \right) \leq \left( \kappa
_{0}+\varepsilon \right) d_{H}\left( U_{1},U_{2}\right) .  \label{eq_14}
\end{equation}%
Once this is proved, we will have $\mathrm{lip\,}\mathcal{F}\left(
U_{0},x_{0}\right) \leq \mathrm{lip\,}\mathcal{F}^{T}\left( \sigma
_{0},x_{0}\right) .$ Let $U_{1},$ $U_{2},$ and $x^{1}$ be given as above.
Associated with the pair $\left( U_{1},U_{2}\right) $ consider the pair of
indexations $\left( \sigma _{1},\sigma _{2}\right) \in \left( \mathbb{R}%
^{n+1}\right) ^{T}\times \left( \mathbb{R}^{n+1}\right) ^{T}$ as in the
previous lemma. Then, we have $d_{\infty }\left( \sigma _{1},\sigma
_{2}\right) =d_{H}\left( U_{1},U_{2}\right) $ and%
\begin{equation*}
d_{\infty }\left( \sigma _{i},\sigma _{0}\right) \leq 3\max \{d_{H}\left(
U_{1},U_{0}\right) ,d_{H}\left( U_{2},U_{0}\right) <\delta .
\end{equation*}

Moreover, it is clear that $\mathrm{rge\,}\sigma _{i}=U_{i},$ for $i=1,2.$
More in detail, the inclusion $\subset $ is evident since we are selecting
projections on $U_{i},$ and $\supset $ comes from $\sigma _{i}\left(
t\right) =t$ for $t\in U_{i}.$ Then, the aimed result will follow
straightforwardly from (\ref{eq_12}). Specifically,%
\begin{equation*}
d\left( x^{1},\mathcal{F}\left( U_{2}\right) \right) =d\left( x^{1},\mathcal{%
F}^{T}\left( \sigma _{2}\right) \right) \leq \left( \kappa _{0}+\varepsilon
\right) d_{\infty }\left( \sigma _{1},\sigma _{2}\right) =\left( \kappa
_{0}+\varepsilon \right) d_{H}\left( U_{1},U_{2}\right) .
\end{equation*}%
This finishes the proof of $\mathrm{lip\,}\mathcal{F}\left(
U_{0},x_{0}\right) \leq \mathrm{lip\,}\mathcal{F}^{T}\left( \sigma
_{0},x_{0}\right) .$

The opposite inequality follows from (\ref{q_0006}). More in detail, assume
that (\ref{eq_14}) holds for all $U_{1},U_{2}\in CL\left( \mathbb{R}%
^{n+1}\right) $ with $d_{\infty }\left( U_{i},U_{0}\right) <\delta ,$ $%
i=1,2, $ and all $x^{1}\in \mathcal{F}\left( U_{1}\right) $ with $\left\Vert
x^{1}-x_{0}\right\Vert <\delta ,$ for some $\varepsilon >0$ and some
associated $\delta >0$; and, for the same $\delta $, consider any pair $%
\sigma _{1},\sigma _{2}\in \left( \mathbb{R}^{n+1}\right) ^{T}$ with $%
d_{\infty }\left( \sigma _{i},\sigma _{0}\right) <\delta ,$ $i=1,2,$ and any 
$x^{1}\in \mathcal{F}^{T}\left( \sigma _{1}\right) .$ Then, appealing to (%
\ref{q_0006}), we conclude from (\ref{eq_14}) that 
\begin{eqnarray*}
d\left( x^{1},\mathcal{F}^{T}\left( \sigma _{2}\right) \right) &=&d\left(
x^{1},\mathcal{F}\left( \mathrm{cl\,rge\,}\sigma _{2}\right) \right) \\
&\leq &\left( \kappa _{0}+\varepsilon \right) d_{H}\left( \mathrm{cl\,rge\,}%
\sigma _{1},\mathrm{cl\,rge\,}\sigma _{2}\right) \\
&\leq &\left( \kappa _{0}+\varepsilon \right) d_{\infty }\left( \sigma
_{1},\sigma _{2}\right) .
\end{eqnarray*}%
\bigskip
\end{dem}

The following lemma constitutes the counterpart of \cite[Thm. 1, Lem. 2]%
{CGHP19} in the context of systems (\ref{eq_problemRO})$.$ We omit the proof
since it follows straightforwardly from the original reference (for systems (%
\ref{eq_system_T})), as far as it only involves a fix system (it is not of
parametric nature). We say that $U\in CL\left( \mathbb{R}^{n+1}\right) $
satisfies the strong Slater condition (SSC, in brief) when $\widehat{x}\in 
\mathbb{R}^{n}$ exists such that $\sup_{\left( a,b\right) \in U}\left(
a^{\prime }\widehat{x}-b\right) <0;$ in such a case $\widehat{x}$ is called
a strong Slater point of $U.$

\begin{lem}
\label{Lem_SSC}We have that

$\left( i\right) $ $U_{0}$ satisfies the SSC if and only if $0_{n}\notin
C_{U_{0}};$

$\left( ii\right) $ Assume that $\left\{ a\in \mathbb{R}^{n}:\left(
a,b\right) \in U_{0}\right\} $ is bounded. Then, $x_{0}$ is a strong Slater
point of $U_{0}$ if and only if $C_{U_{0}}=\emptyset .$
\end{lem}

As a consequence of the previous lemma, we derive the following corollary
which gather two special particular cases.

\begin{rem}
\emph{Note that under the assumptions of Theorem \ref{theo_lip_U}:}

$\left( i\right) $ $\mathrm{lip\,}\mathcal{F}\left( U_{0},x_{0}\right)
=+\infty $ \emph{if and only if} $0_{n}\in C_{U_{0}},$ \emph{equivalently
SSC is not satisfied at} $U_{0}$ \emph{.}

$\left( ii\right) $ $\mathrm{lip\,}\mathcal{F}\left( U_{0},x_{0}\right) =0$ 
\emph{if and only if} $C_{U_{0}}=\emptyset ,$ \emph{equivalently, }$x_{0}$ 
\emph{is an SS element}$.$
\end{rem}

\begin{cor}
\label{Corol_1}\emph{(see \cite[Proposition 1]{CGP08}) }Let $U_{0}\in
CL\left( \mathbb{R}^{n+1}\right) $ be such that $\left\{ a\in \mathbb{R}%
^{n}:\left( a,b\right) \in U_{0}\right\} $ is bounded and assume that SSC
holds at $U_{0}.$ Then $\mathcal{F}$ has the Aubin property at $\left(
U_{0},x_{0}\right) $ for any $x_{0}\in \mathcal{F}\left( U_{0}\right) $.
\end{cor}

\section{Application to convex inequalities}

This section is devoted to apply the previous results about linear systems
to the convex case. Throughout this section $\mathbb{R}^{n}$ is considered
to be endowed with the Euclidean norm, denoted in the same way for
simplicity, $\left\Vert \cdot \right\Vert ,$ and $\mathbb{B}$ is the
corresponding closed unit ball.

We consider the parameterized family of convex inequalities (\ref%
{eq_convex_ineq}) and the corresponding feasible set mapping $\mathcal{L}%
:\Gamma \rightrightarrows \mathbb{R}^{n}$ assigning to each convex function $%
f\in \Gamma $ its zero (sub)level set%
\begin{equation*}
\mathcal{L}\left( f\right) :=\left\{ x\in \mathbb{R}^{n}:f\left( x\right)
\leq 0\right\} ,\text{ }f\in \Gamma .
\end{equation*}%
It is well-known that, for each $f\in \Gamma ,$ $\mathcal{L}\left( f\right)
\subset \mathbb{R}^{n}$ is a closed convex set and, as commented in Section
1, via a standard linearization, it can be written as the feasible set of a
linear semi-infinite inequality system of the form; i.e., 
\begin{equation}
\mathcal{L}\left( f\right) =\left\{ x\in \mathbb{R}^{n}:a^{\prime }x\leq
a^{\prime }z-f\left( z\right) ,\text{~~}\left( z,a\right) \in \mathrm{gph}%
\partial f\right\} .  \label{eq_004}
\end{equation}%
First, let us see that we can reduce the index set of system (\ref{eq_004})
to a certain subset of $\mathrm{gph}\partial f$.

\begin{lem}
\label{lem_basic1}Let $f\in \Gamma $ and $X\subset \mathbb{R}^{n}$ be an
open set such that $\mathcal{L}\left( f\right) \subset X.$ Then,%
\begin{equation}
\mathcal{L}\left( f\right) =\left\{ x\in \mathbb{R}^{n}:a^{\prime }x\leq
a^{\prime }z-f\left( z\right) ,\text{~~}a\in \partial f\left( z\right) ,%
\text{ }z\in X\right\} .  \label{eq_005}
\end{equation}

Moreover, $X$ can be replaced in (\ref{eq_005}) with any $S\supset X$.
\end{lem}

\begin{dem}
The inclusion `$\subset $' is trivial. Let us prove `$\supset $' reasoning
by contradiction. Assume the existence of $x_{0}\in \mathbb{R}^{n}$ such
that 
\begin{equation}
a^{\prime }x_{0}\leq a^{\prime }z-f\left( z\right) ,~\text{whenever }a\in
\partial f\left( z\right) ,\text{ with }z\in X,  \label{eq_006}
\end{equation}%
and $x_{0}\notin \mathcal{L}\left( f\right) ,$ which entails $f\left(
x_{0}\right) >0.$ Indeed, we have $x_{0}\notin X;$ otherwise, taking $%
z=x_{0} $ in (\ref{eq_006}), we would have $a%
{\acute{}}%
x_{0}\leq a%
{\acute{}}%
x_{0}-f\left( x_{0}\right) ,$ for any $a\in \partial f\left( x_{0}\right) $,
yielding the contradiction $f\left( x_{0}\right) \leq 0.$ Once we know that $%
x_{0}\notin X,$ pick any $x_{1}\in \mathcal{L}\left( f\right) \subset X$ and
define $x^{\lambda }:=\left( 1-\lambda \right) x_{1}+\lambda x_{0}$ for any $%
0<\lambda <1.$ Observe that for each $0<\lambda <1,$ $x^{\lambda }$ also
verifies the linear inequalities of the right member in (\ref{eq_006}), i.
e., 
\begin{equation*}
a^{\prime }x^{\lambda }\leq a^{\prime }z-f\left( z\right) ,~\text{whenever }%
a\in \partial f\left( z\right) ,\text{ with }z\in X.
\end{equation*}%
Then, arguing as in the previous paragraph, $x^{\lambda }\notin X\,,$ $%
0<\lambda <1$, which represents a contradiction since we can choose $%
x^{\lambda }$ sufficiently close to $x^{1}$ to ensure $x^{\lambda }\in X.$

Finally, from (\ref{eq_004}) and (\ref{eq_005}), it is obvious that, $X$ can
be replaced by any $S\supset X$.
\end{dem}

From now on we use the notation: $f_{0}\in \Gamma $ is our nominal convex
function$,$ $\alpha _{0}>0$ is a fixed scalar, and $E_{0}\subset \mathbb{R}%
^{n}$ is the $\alpha _{0}$-enlargement of the nominal feasible set $\mathcal{%
L(}f_{0});$ i.e., 
\begin{equation}
E_{0}:=\mathcal{L(}f_{0})+\alpha _{0}\mathbb{B}.  \label{eq_0010}
\end{equation}%
As a particular consequence of the previous lemma, we can write%
\begin{equation*}
\mathcal{L}\left( f_{0}\right) =\left\{ x\in \mathbb{R}^{n}:a^{\prime }x\leq
a^{\prime }z-f_{0}\left( z\right) ,\text{~~}a\in \partial f_{0}\left(
z\right) ,\text{ }z\in E_{0}\right\} .
\end{equation*}%
Going further, the following lemma ensures that we can keep the same $E_{0}$
in the linear representation of $\mathcal{L(}f_{1})$ provided that $f_{1}\in
\Gamma $ is close enough to $f_{0}$ in relation to the pseudo-distance $%
d_{E_{0}}$ defined in (\ref{eq_0007}).

\begin{lem}
\label{lem_basic2}Assume that $\mathcal{L(}f_{0})$ is bounded$.$ There
exists $\eta >0$ such that 
\begin{equation*}
\mathcal{L(}f_{1})=\{x\in \mathbb{R}^{n}:u^{\prime }x\leq u^{\prime
}z-f_{1}\left( z\right) ,\text{~}u\in \partial f_{1}\left( z\right) ,\text{ }%
z\in E_{0}\mathbb{\}},
\end{equation*}%
whenever $f_{1}\in \Gamma ,$ with $d_{E_{0}}\left( f_{1},f_{0}\right) <\eta
. $

($E_{0}$ can be replaced by any set $E\supset E_{0}$).
\end{lem}

\begin{dem}
Obviously $\mathcal{L}\left( f_{0}\right) +\frac{\mathcal{\alpha }_{0}}{2}%
\mathbb{B}$ is a compact convex subset of $\mathbb{R}^{n}$ and so the
following minimum is attained:%
\begin{equation*}
m:=\min \left\{ f_{0}\left( x\right) :x\in \mathrm{bd}\,\left( \mathcal{L}%
\left( f_{0}\right) +\frac{\mathcal{\alpha }_{0}}{2}\mathbb{B}\right)
\right\} =f_{0}\left( \widetilde{x}\right) ,
\end{equation*}%
for some $\widetilde{x}\in \mathrm{bd}\,\left( \mathcal{L}\left(
f_{0}\right) +\frac{\mathcal{\alpha }_{0}}{2}\mathbb{B}\right) .$ Observe
that $m>0,$ since otherwise we would have $\widetilde{x}\in \mathcal{L}%
\left( f_{0}\right) ,$ and then $\widetilde{x}\in \widetilde{x}+\frac{%
\mathcal{\alpha }_{0}}{2}\mathbb{B}\subset \mathcal{L}\left( f_{0}\right) +%
\frac{\mathcal{\alpha }_{0}}{2}\mathbb{B},$ which would yield the
contradiction $\widetilde{x}\in \mathrm{int}\,\left( \mathcal{L}\left(
f_{0}\right) +\frac{\mathcal{\alpha }_{0}}{2}\mathbb{B}\right) .$ Take 
\begin{equation*}
\eta :=\frac{m}{4}>0,
\end{equation*}%
and consider any convex function $f_{1}$ such that $d_{E_{0}}\left(
f_{1},f_{0}\right) <\eta .$ Let us see that 
\begin{equation*}
\left\{ x\in \mathbb{R}^{n}:f_{1}\left( x\right) \leq 0\right\} \subset 
\mathcal{L}\left( f_{0}\right) +\frac{\mathcal{\alpha }_{0}}{2}\mathbb{B}.
\end{equation*}

For any $x\in \mathrm{bd}\,\left( \mathcal{L}\left( f_{0}\right) +\frac{%
\mathcal{\alpha }_{0}}{2}\mathbb{B}\right) $ we have%
\begin{equation}
f_{1}\left( x\right) =f_{0}\left( x\right) +f_{1}\left( x\right)
-f_{0}\left( x\right) >m-\frac{m}{4}=\frac{3m}{4},  \label{eq_21001}
\end{equation}%
while, for $x\in \mathcal{L}\left( f_{0}\right) $ one has%
\begin{equation*}
f_{1}\left( x\right) \leq f_{1}\left( x\right) -f_{0}\left( x\right) <\frac{m%
}{4}.
\end{equation*}

Now, arguing by contradiction, assume that there exists $x_{1}\in \mathbb{R}%
^{n}$ such that $f_{1}\left( x_{1}\right) \leq 0$ and $x_{1}\notin \mathcal{L%
}\left( f_{0}\right) +\frac{\mathcal{\alpha }_{0}}{2}\mathbb{B}.$ Take any $%
x_{0}\in \mathcal{L}\left( f_{0}\right) $ (in particular, $x_{0}\in \mathrm{%
int}\,\left( \mathcal{L}\left( f_{0}\right) +\frac{\mathcal{\alpha }_{0}}{2}%
\mathbb{B}\right) $) and let $\lambda \in \left] 0,1\right[ $ such that 
\begin{equation*}
\left( 1-\lambda \right) x_{0}+\lambda x_{1}\in \mathrm{bd}\,\left( \mathcal{%
L}\left( f_{0}\right) +\frac{\mathcal{\alpha }_{0}}{2}\mathbb{B}\right) .
\end{equation*}%
Then, we attain the contradiction, with (\ref{eq_21001}), 
\begin{equation*}
f_{1}\left( \left( 1-\lambda \right) x_{0}+\lambda x_{1}\right) \leq \left(
1-\lambda \right) f_{1}\left( x_{0}\right) +\lambda f_{1}\left( x_{1}\right)
<\left( 1-\lambda \right) \frac{m}{4}<\frac{m}{4}.
\end{equation*}%
The fact that $E_{0}$ can be replaced by any subset containing it comes from
the standard linearization of the convex inequality $f_{1}\left( x\right)
\leq 0$ where the whole graph, $\mathrm{gph}\partial f_{1},$ is used.
\end{dem}

The following lemma constitutes a key tool for our purposes.

\begin{lem}
\label{Lem3}Let $K_{1},K_{2}\subset \mathbb{R}^{n}$ be compact sets, let $%
f_{1},f_{2}:\mathbb{R}^{n}\rightarrow \mathbb{R}$ be convex functions, and
consider%
\begin{equation*}
U_{i}:=\left\{ \left( a,b\right) :b=a^{\prime }z-f_{i}\left( z\right) ,\text{
}a\in \partial f_{i}\left( z\right) ,\text{ }z\in K_{i}\right\} .\,\ i=1,2.
\end{equation*}%
Then,%
\begin{equation}
d_{H}\left( U_{1},U_{2}\right) \leq \rho \text{ }d_{H}\left( \partial
f_{1}\left( K_{1}\right) ,\partial f_{2}\left( K_{2}\right) \right)
+d_{K_{1}\cup K_{2}}\left( f_{1},f_{2}\right) ,  \label{eq_2000}
\end{equation}%
where $\rho :=\max \{1+\left\Vert x\right\Vert :x\in K_{1}\cup K_{2}\}.$
\end{lem}

\begin{dem}
Take $K_{i},f_{i},U_{i},$ for $i=1,2,$ and $\rho ,$ as in the statement of
the lemma$.$ Let us establish the inequality 
\begin{equation}
e\left( U_{1},U_{2}\right) \leq \rho \text{ }d_{H}\left( \partial
f_{1}\left( K_{1}\right) ,\partial f_{2}\left( K_{2}\right) \right)
+d_{K_{1}\cup K_{2}}\left( f_{1},f_{2}\right) ,  \label{eq_50}
\end{equation}%
which yields by symmetry the aimed inequality (\ref{eq_2000}). For
simplicity, in this proof we use the notation 
\begin{equation*}
\xi :=d_{H}\left( \partial f_{1}\left( K_{1}\right) ,\partial f_{2}\left(
K_{1}\right) \right) .
\end{equation*}%
Specifically, take any $\left( a_{1},b_{1}\right) \in U_{1},$ and let us
prove the existence of $\left( a_{2},b_{2}\right) \in U_{2}$ such that 
\begin{equation}
\left\Vert \left( a_{1},b_{1}\right) -\left( a_{2},b_{2}\right) \right\Vert
\leq \rho \text{ }\xi +d_{K_{1}\cup K_{2}}\left( f_{1},f_{2}\right) .
\label{eq_3000}
\end{equation}

By definition, $\left( a_{1},b_{1}\right) \in U_{1}$ entails the existence
of $z_{1}\in K_{1}$ such that%
\begin{equation*}
a_{1}\in \partial f_{1}\left( z_{1}\right) \text{ and }b_{1}=a_{1}^{\prime
}z_{1}-f_{1}\left( z_{1}\right) .
\end{equation*}%
Since $\partial f_{1}\left( K_{1}\right) $ and $\partial f_{2}\left(
K_{1}\right) $ are compact subsets in $\mathbb{R}^{n}$ (see again \cite[%
Theorem 24.7]{Rock}), in particular we have that 
\begin{equation}
\partial f_{1}\left( K_{1}\right) \subset \partial f_{2}\left( K_{2}\right)
+\xi \mathbb{B},  \label{eq_30}
\end{equation}%
and, so, we can write%
\begin{equation*}
a_{1}=a_{2}+\xi w,\text{ for some }a_{2}\in \partial f_{2}\left(
z_{2}\right) ,\text{ }z_{2}\in K_{2},\text{ and }\left\Vert w\right\Vert
\leq 1.
\end{equation*}

Define 
\begin{equation*}
b_{2}:=a_{2}^{\prime }z_{2}-f_{2}\left( z_{2}\right) ,
\end{equation*}%
and let us establish (\ref{eq_3000}) for such an element $\left(
a_{2},b_{2}\right) $.

On the one hand, 
\begin{eqnarray*}
a_{1}^{\prime }z_{1}-f_{1}\left( z_{1}\right) &=&a_{2}^{\prime }z_{1}+\xi
w^{\prime }z_{1}-f_{1}\left( z_{1}\right) \\
&=&a_{2}^{\prime }z_{1}-f_{2}\left( z_{1}\right) +\xi w^{\prime
}z_{1}+f_{2}\left( z_{1}\right) -f_{1}\left( z_{1}\right) \\
&\leq &a_{2}^{\prime }z_{2}-f_{2}\left( z_{2}\right) +\xi \left\Vert
z_{1}\right\Vert +f_{2}\left( z_{1}\right) -f_{1}\left( z_{1}\right) \\
&\leq &a_{2}^{\prime }z_{2}-f_{2}\left( z_{2}\right) +\xi \rho +d_{K_{1}\cup
K_{2}}\left( f_{1},f_{2}\right) ,
\end{eqnarray*}%
where for the first inequality we have applied the fact that $a_{2}\in
\partial f_{2}\left( z_{2}\right) .$

On the other hand, since $a_{1}\in \partial f_{1}\left( z_{1}\right) ,$ we
have 
\begin{eqnarray*}
a_{1}^{\prime }z_{1}-f_{1}\left( z_{1}\right) &\geq &a_{1}^{\prime
}z_{2}-f_{1}\left( z_{2}\right) \\
&=&a_{2}^{\prime }z_{2}+\xi w^{\prime }z_{2}-f_{1}\left( z_{2}\right) \\
&=&a_{2}^{\prime }z_{2}-f_{2}\left( z_{2}\right) +\xi w^{\prime
}z_{2}+f_{2}\left( z_{2}\right) -f_{1}\left( z\right) \\
&\geq &a_{2}^{\prime }z_{2}-f_{2}\left( z_{2}\right) -\xi \rho -d_{K_{1}\cup
K_{2}}\left( f_{1},f_{2}\right) .
\end{eqnarray*}

So, we have established%
\begin{equation*}
\left\vert b_{1}-b_{2}\right\vert \leq \xi \rho +d_{K_{1}\cup K_{2}}\left(
f_{1},f_{2}\right) .
\end{equation*}%
Finally, we have (recall that $\rho >1),$ 
\begin{eqnarray*}
\left\Vert \left( a_{1},b_{1}\right) -\left( a_{2},b_{2}\right) \right\Vert
&=&\max \{\left\Vert a_{1}-a_{2}\right\Vert ,\left\vert
b_{1}-b_{2}\right\vert \} \\
&\leq &\max \{\xi ,\rho \xi +d_{K_{1}\cup K_{2}}\left( f_{1},f_{2}\right) \\
&=&\rho \xi +d_{K_{1}\cup K_{2}}\left( f_{1},f_{2}\right) ,
\end{eqnarray*}%
which yields (\ref{eq_3000}) and the proof is complete.
\end{dem}

The following theorem provides the announced results about the Lipschitzian
behavior of the feasible set of convex inequalities. It appeals to the
constant%
\begin{equation}
\kappa _{0}:=\frac{\left\Vert x_{0}\right\Vert +1}{d\left(
0_{n},C_{U_{0}}\right) },  \label{eq_009}
\end{equation}%
where 
\begin{equation*}
U_{0}:=\left\{ \left( a,b\right) :b=a%
{\acute{}}%
z-f_{0}\left( z\right) ,\text{ }a\in \partial f_{0}\left( z\right) ,\text{ }%
z\in E_{0}\right\} .
\end{equation*}%
Before the theorem, the next proposition says that $\kappa _{0}$ is finite
when the convex inequality `$f_{0}\left( x\right) \leq 0$' has a strict
solution (i.e., when this convex inequality verifies the Slater condition).

\begin{prop}
There exists $z_{0}\in \mathbb{R}^{n}$ such that $f_{0}\left( z_{0}\right)
<0 $ if and only if $0_{n}\notin C_{U_{0}}.$
\end{prop}

\begin{dem}
According to Lemma \ref{Lem_SSC}, we only have to prove that the existence
of $z_{0}\in \mathbb{R}^{n}$ such that $f_{0}\left( z_{0}\right) <0$ is
equivalent to SSC at $U_{0}.$

Take $z_{0}\in \mathbb{R}^{n}$ such that $f_{0}\left( z_{0}\right) <0$ and
let us see that $z_{0}$ is a SS point of $U_{0}.$ Observe that%
\begin{equation*}
0>f_{0}\left( z_{0}\right) \geq f_{0}\left( z\right) +a^{\prime }\left(
z_{0}-z\right) ,\text{ whenever }a\in \partial f_{0}\left( z\right) ,\text{ }%
z\in E_{0};
\end{equation*}%
equivalently%
\begin{equation*}
a^{\prime }z_{0}\leq a^{\prime }z-f_{0}\left( z\right) +f_{0}\left(
z_{0}\right) ,\text{ for all }a\in \partial f_{0}\left( z\right) ,\text{
with }z\in E_{0}.
\end{equation*}%
So,\ $\sup_{\binom{a}{v}\in U_{0}}\left( a^{\prime }z_{0}-v\right) \leq
f_{0}\left( z_{0}\right) <0.$

Reciprocally, let $z_{0}$ be a SS point of $U_{0}$, in particular, $z_{0}\in 
\mathcal{L(}f_{0})\subset E_{0}$ and, so, taking any $a\in \partial
f_{0}\left( z_{0}\right) ,$ we have that 
\begin{equation*}
a^{\prime }z_{0}<a^{\prime }z_{0}-f_{0}\left( z_{0}\right) ,
\end{equation*}%
which entails $f_{0}\left( z_{0}\right) <0.$
\end{dem}

\begin{theo}
\label{Th Lipschitz convex}Let $\kappa >\kappa _{0},$ $\alpha >0,$ $%
E:=E_{0}+\alpha \mathbb{B},$ and $\rho :=\max \{1+\left\Vert x\right\Vert
:x\in E\}.$ Then, there exist $\delta _{0}>0$ such that $0<\delta \leq
\delta _{0}$ implies%
\begin{equation*}
d\left( x_{1},\mathcal{L(}f_{2})\right) \leq \kappa \left( \rho d_{H}\left(
\partial f_{1}\left( E_{0}+\sqrt{\delta }\mathbb{B}\right) ,\partial
f_{2}\left( E_{0}+\sqrt{\delta }\mathbb{B}\right) \right) +d_{E}\left(
f_{1},f_{2}\right) \right) ,\text{ }
\end{equation*}%
provided that $f_{1},f_{2}\in \Gamma $, $x_{1}\in \mathcal{L(}f_{1}),$ $%
d_{_{E}}\left( f_{i},f_{0}\right) \leq \delta ,$ $i=1,2,$ and $\left\Vert
x_{1}-x_{0}\right\Vert \leq \delta .$
\end{theo}

\begin{dem}
Take $\kappa >\kappa _{0}$ and fix $\varepsilon >0.$ Theorem \ref{theo_lip_U}
ensures the existence of $\delta _{1}>0$ such that $d_{H}\left(
U_{i},U\right) \leq \delta _{1},$ $U_{i}\in CL\left( \mathbb{R}^{n+1}\right)
,$ $i=1,2,...,$ and $\left\Vert x_{1}-x_{0}\right\Vert \leq \delta ,$ $%
x_{1}\in \mathcal{F}\left( U_{1}\right) ,$ imply 
\begin{equation}
d\left( x_{1},\mathcal{F}\left( U_{2}\right) \right) \leq \kappa d_{H}\left(
U_{1},U_{2}\right) .  \label{eq_0001}
\end{equation}%
On the other hand, according to Corollary \ref{Cor_subdifferential}$\left(
iii\right) ,$ choose $\delta _{2}>0$ such that $d_{_{E}}\left(
f,f_{0}\right) \leq \delta \leq \delta _{2}$ implies $,$%
\begin{equation*}
d_{H}\left( \partial f\left( E_{0}+\sqrt{\delta }\mathbb{B}\right) ,\partial
f_{0}\left( E_{0}\right) \right) \leq \dfrac{\delta _{1}}{2\rho }.
\end{equation*}%
Let $\eta >0$ be as in Lemma \ref{lem_basic2}, and consider%
\begin{equation*}
\delta _{0}=\min \{\dfrac{\delta _{1}}{2},\delta _{2},\eta \}.
\end{equation*}%
Now consider $0<\delta \leq \delta _{0},$ and $f_{1},f_{2}\in \Gamma $, with 
$d_{_{E}}\left( f_{i},f_{0}\right) \leq \delta .$ Define, 
\begin{equation*}
U_{i}^{\delta }:=\left\{ \left( a,b\right) :b=a^{\prime }z-f_{i}\left(
z\right) ,\text{ }a\in \partial f_{i}\left( z\right) ,\text{ }z\in E_{0}+%
\sqrt{\delta }\mathbb{B}\right\} \,\ i=1,2.
\end{equation*}%
Appealing to Lemma \ref{Lem3}, we have, for $i=1,2,$%
\begin{eqnarray*}
d_{H}\left( U_{i}^{\delta },U_{0}\right) &\leq &\rho d_{H}\left( \partial
f_{i}\left( E_{0}+\sqrt{\delta }\mathbb{B}\right) ,\partial f_{0}\left(
E_{0}\right) \right) +d_{E_{0}+\sqrt{\delta }\mathbb{B}}\left(
f_{1},f_{2}\right) \\
&\leq &\rho \dfrac{\delta _{1}}{2\rho }+d_{E}\left( f_{1},f_{2}\right) \leq 
\dfrac{\delta _{1}}{2}+\dfrac{\delta _{1}}{2}=\delta _{1}.
\end{eqnarray*}%
Moreover, since $d_{_{E}}\left( f,f_{0}\right) \leq \delta \leq \eta ,$ we
have%
\begin{equation*}
\mathcal{F}\left( U_{i}^{\delta }\right) =\mathcal{L}\left( f_{i}\right) ,%
\text{ }i=1,2.
\end{equation*}

Consequently, appealing to (\ref{eq_0001}) in the particular case $%
U_{i}=U_{i}^{\delta },$ $i=1,2,$ we conclude%
\begin{eqnarray*}
d\left( x_{1},\mathcal{L(}f_{2})\right) &\leq &\kappa d_{H}\left(
U_{1}^{\delta },U_{2}^{\delta }\right) \\
&\leq &\kappa \left( \rho d_{H}\left( \partial f_{1}\left( E_{0}+\sqrt{%
\delta }\mathbb{B}\right) ,\partial f_{2}\left( E_{0}+\sqrt{\delta }\mathbb{B%
}\right) \right) +d_{E}\left( f_{1},f_{2}\right) \right) ,
\end{eqnarray*}%
where in the second inequality we have appealed again to Lemma \ref{Lem3}.
\end{dem}

\subsection{The convex differentiable case}

Throughout this subsection we assume that our nominal function $f_{0}\in
\Gamma $ is differentiable, so that we write $\nabla f_{0}$ instead of $%
\partial f_{0}$. The following theorem provides the counterpart of Corollary %
\ref{Cor_subdifferential}$\left( iii\right) $ under differentiability of $%
f_{0}$.

\begin{theo}
Let $K_{0}\subset \mathbb{R}^{n}$ a compact set, $\alpha >0,$ and $%
K:=K_{0}+\alpha \mathbb{B}.$ Given $\varepsilon >0,$ there exists $\delta >0$
such that for any $f\in \Gamma ,$ with $d_{K}\left( f,f_{0}\right) \leq
\delta ,$ one has%
\begin{equation*}
d_{H}\left( \partial f\left( K_{0}\right) ,\nabla f_{0}\left( K_{0}\right)
\right) \leq \varepsilon .
\end{equation*}
\end{theo}

\begin{dem}
Take $\varepsilon >0.$ From Theorem \ref{Th_Subdif} $\left( i\right) \ $%
there exists $\delta _{1}>0$ such that 
\begin{equation*}
\partial f\left( K_{0}+\delta _{1}\mathbb{B}\right) \subset \nabla
f_{0}\left( K_{0}\right) +\varepsilon \mathbb{B},\text{ provided that }%
d_{K}\left( f,f_{0}\right) \leq \delta _{1},\text{ }f\in \Gamma .
\end{equation*}%
In particular, $\partial f\left( K_{0}\right) \subset \nabla f_{0}\left(
K_{0}\right) +\varepsilon \mathbb{B},$ if $d_{K}\left( f,f_{0}\right) \leq
\delta _{1},$ $f_{1}\in \Gamma .$

Let us prove the existence of $\delta _{2}>0\,$such that 
\begin{equation*}
\nabla f_{0}\left( K_{0}\right) \subset \partial f\left( K_{0}\right)
+\varepsilon \mathbb{B},\text{ if }d_{K}\left( f,f_{0}\right) \leq \delta
_{2},\text{ }f\in \Gamma .
\end{equation*}%
In such a case, just take $\delta :=\min \{\delta _{1},\delta _{2}\}$ to
finish the proof.

Arguing by contradiction, assume the existence of a sequence of convex
functions $\{f_{r}\},$ with $d_{K}\left( f_{r},f_{0}\right) \leq \frac{1}{r}$
such that 
\begin{equation*}
\nabla f_{0}\left( K_{0}\right) \not\subset \partial f_{r}\left(
K_{0}\right) +\varepsilon \mathbb{B},\text{ for all }r.
\end{equation*}%
For each $r,$ let $x^{r}\in K_{0}$ such that 
\begin{equation}
\nabla f_{0}\left( x^{r}\right) \notin \partial f_{r}\left( K_{0}\right)
+\varepsilon \mathbb{B}.  \label{eq_0002}
\end{equation}%
The compactness of $K_{0},$ and consequently of $\nabla f_{0}\left(
K_{0}\right) $ allows us to assume that $\{x^{r}\}$ and $\{\nabla
f_{0}\left( x^{r}\right) \}\,$converge to $\overline{x}\in K_{0}$ and $%
\nabla f_{0}\left( \overline{x}\right) ,$ respectively (see again \cite[%
Theorem 24.5]{Rock}). This fact, together with (\ref{eq_0002}) yields the
existence of $r_{0}\in \mathbb{N}$ such that%
\begin{equation}
\left( \nabla f_{0}\left( \overline{x}\right) +\frac{\varepsilon }{2}\mathbb{%
B}\right) \cap \partial f_{r}\left( K_{0}\right) =\emptyset ,\text{ for }%
r\geq r_{0}.  \label{eq_002}
\end{equation}%
On the other hand, \cite[Theorem 24.5]{Rock} guarantees, for $r$ large
enough, 
\begin{equation*}
\partial f_{r}\left( \overline{x}\right) \subset \nabla f_{0}\left( 
\overline{x}\right) +\frac{\varepsilon }{2}\mathbb{B},
\end{equation*}%
which represents a contradiction.
\end{dem}

Following the proof of Theorem \ref{Th Lipschitz convex}, appealing to the
previous theorem instead of Corollary \ref{Cor_subdifferential}$\left(
iii\right) $ we derive the following corollary. Recall that $E_{0}$ and $%
k_{0}$ are defined in (\ref{eq_0010}) and (\ref{eq_009}), respectively.

\begin{cor}
Let $\kappa >\kappa _{0},$ $\alpha >0,$ $E:=E_{0}+\alpha \mathbb{B},$ and $%
\rho :=\max \{1+\left\Vert x\right\Vert :x\in E\}.$ There exist $\delta >0$
such that for any $f_{1},f_{2}\in \Gamma $, with $d_{_{E}}\left(
f_{i},f_{0}\right) \leq \delta ,$ $i=1,2,$ and any $x^{1}\in \mathcal{L(}%
f_{1}),$ with $\left\Vert x^{1}-x_{0}\right\Vert \leq \delta ,$ one has 
\begin{equation}
d\left( x^{1},\mathcal{L(}f_{2})\right) \leq \kappa \left( \rho d_{H}\left(
\partial f_{1}\left( E_{0}\right) ,\partial f_{2}\left( E_{0}\right) \right)
+d_{E}\left( f_{1},f_{2}\right) \right) .  \label{eq_0008}
\end{equation}
\end{cor}

\end{document}